\documentclass[11 pt]{amsart}
\usepackage[T1]{fontenc}
\usepackage[utf8]{inputenc}
\usepackage[english]{babel}
\usepackage{mathtools}
\usepackage{amsfonts}
\usepackage{amsthm}
\usepackage{braket}
\usepackage{lipsum}
\numberwithin{equation}{section}
\theoremstyle{plain}
\newtheorem*{theo*}{Theorem}
\newtheorem{theo}{Theorem}[section]
\newtheorem{coro}[theo]{Corollary}
\theoremstyle{definition}
\newtheorem{ex}[theo]{Example}
\newtheorem{obs}[theo]{Observation}
\newtheorem{defi}[theo]{Definition}

\newtheorem{lemma}[theo]{Lemma}
\newcommand{\N}{\mathbb{N}}
\newcommand{\C}{\mathbb{C}}
\newcommand{\D}{\mathbb{D}}
\newcommand{\de}{\partial}
\newcommand{\cp}{\mathrm{cap}}
\newcommand{\s}{\mathcal{S}}

\DeclarePairedDelimiter{\abs}{\lvert}{\rvert}
\DeclarePairedDelimiter{\babs}{\bigg\lvert}{\bigg\rvert}

\DeclarePairedDelimiter\floor{\lfloor}{\rfloor}
\title{Zeros of normalized sections of non convergent power series}
\author{Alberto Dayan}
\address{Department of Mathematics\newline Washington University in St. Louis,\newline One Brookings Drive, St. Louis, MO 63130, USA}
\email{alberto.dayan@wustl.edu}
\date{\today}
\begin{document}

\maketitle
\begin{abstract}
A well known result due to Carlson \cite{carlo} affirms that a power series with finite and positive radius of convergence $R$ has no Ostrowski gaps if and only if the sequence of zeros of its $n$th sections is asymptotically equidistributed to $\de\D_R$. Here we extend this characterization to those power series with null radius of convergence, modulo some necessary normalizations of the sequence of the sections of $f$.
\end{abstract}
\section{Introduction}
\protect\footnote{Research supported by Starting Grant (StG), PE1, ERC-2012-StG-20111012}Let
\begin{equation}
\label{eqn:ps}
f(z)=\sum_{k=0}^\infty a_kz^k
\end{equation}
be a complex-valued power series with $n$-th partial sum
\[
s_n(z)=s_n(f)(z):=\sum_{k=0}^n a_kz^k.
\]
If $Z_n(f)$ is the multi-set of the zeros of $s_n$, define the $n$-th \emph{zero counting measure} of $f$ 
\[
\nu_n=\nu_n(f):=\frac{1}{n}\sum_{w\in Z_n(f)}\delta_w
\]
as the weighted sum of Dirac deltas placed at the zeros of $s_n$, repeated according to their multiplicity. If $a_n=0$, we set $n-\deg(s_n)$ zeros of $s_n$ at $\infty_{\C}$, so that $\nu_n$ is, for any natural $n$, a probability measure on the Riemann sphere $\overline{\C}:=\C\cup\set{\infty_\C}$.\\
The asymptotic distribution of the zeros of the sections $s_n$ has been intensively investigated, at least in the case in which $f$ has positive and finite radius of convergence $R$, by looking at the (eventual) weak limit of the sequence $(\nu_n)_{n\in\N}$. Jentzsch proved in \cite{jen} that any point of the circle $\de\D_R$ is a limit point of the zeros of $s_n$, and Szegö \cite{szego} improved the result by showing that, in fact, there exists a subsequence of $(\nu_{n_k})_{k\in\N}$ that is asymptotically equidistributed on $\de\D_R$:  if $\Lambda_R$ is the normalized Lebesgue measure on $\de\D_R$,  $(\nu_{n_k})_{k\in\N}$ converges weakly at $\Lambda_R$, i.e., there exists the limit
\[
\lim_{k\to\infty}\frac{1}{n_k}\sum_{w\in Z_{n_k}(f)}h(w)=\frac{1}{2\pi R}\int_0^{2\pi}h(Re^{i\theta})\,d\theta
\]
for any function $h$ bounded and continuous on $\overline{\C}$. Finally, Carlson \cite{carlo} proved that the whole sequence $(\nu_n)_{n\in\N}$ converges weakly to $\Lambda_R$ if and only if the power series in \eqref{eqn:ps} has no \emph{Ostrowski gaps}. In this case, we will say that $f$ belongs to the {\bf Szegö class} $\s$. A rigorous definition of Ostrowski gaps is going to be given in Definition \ref{def:ostr}, but essentially the big coefficients of $f$, that is the ones that yield its radius of convergence, appear in any tail of a section $s_n$ of length $\gamma n$, for all $n$ big enough and $\gamma$ in $(0,1)$. \\
As one can notice, instead of taking the partial sums of a predetermined power series, one can consider an arbitrary sequence of polynomials $(p_n)_{n\in\N}$, and look at the behavior of their zero counting measures, even without assuming uniform convergence for $p_n$. In particular, we will say that such a sequence of polynomials belongs to $\s$ if the associated $(\nu_n)_{n\in\N}$ converges weakly to $\Lambda:=\Lambda_1$. As we will see, this will be useful for extending Carlson's result to the case in which $f$ has null radius of convergence. First, we will have to extend the notion of Ostrowski gaps in \cite{erdos} for such a power series: that will be the aim Section \ref{sec:gauge}.  Then we will observe that, in order to hope in an asymptotic equidistribution of the zeros of $f$ to the unit circle, we need to normalize its sections $(s_n)_{n\in\N}$. Finally, we will prove that the sequence of polynomials $(p_n)_{n\in\N}$ obtained is in $\s$ accordingly, again, to a gap condition:
\begin{theo}
\label{theo:main}
Let
\[
f(z)=\sum_{n=0}^\infty a_nz^n
\]
be a power series with null radius of convergence. Then, if
\[
A_n:=\max_{1\leq k\leq n}\abs{a_k}^\frac{1}{k},
\]
then  the sequence of the normalized section of $f$
\begin{equation*}
p_n(z):=s_n\bigg(\frac{z}{A_n} \bigg)=\sum_{k=0}^n a_kA_n^{-k}z^k
\end{equation*}
belongs to the Szegö class if and only if, for any $\gamma$ in $(0, 1)$, 
\[
\liminf_{n\to\infty} \max_{(1-\gamma)n\leq k\leq n}\frac{\abs{a_k}^\frac{1}{k}}{A_n}=1
\]
\end{theo}
 This will extends Carlson's result to non convergent power series. \\
 
 The author thanks José L. Fernández for the great help and support during the development of this work, and John McCarthy for all the conversations about the topic and his valuable suggestions.
\section{Preliminary tools}
This section's aim is to develop the necessary tools for the discussion. First, we will give a notion of \emph{gauge} of a power series with null radius of convergence. After that, we will see some applications of classical bounds regarding zeros of a polynomial and its coefficients. Finally, we will see how potential theory can be applied in order to study the asymptotic distribution of the zero counting measure of a sequence of polynomials $(p_n)_{n\in\N}$.
\subsection{Gauge of sequences}
\label{sec:gauge}
Let $(a_n)_{n\in\N}$ be a sequence of complex numbers, and define, for any $n\geq1$,
\begin{equation}
\label{eqn:alpha}
\alpha_n:=\abs{a_n}^\frac{1}{n}.
\end{equation}
Suppose that 
\[
\limsup_{n\to\infty}\alpha_n=\infty,
\]
i.e., that the sequence is composed by the coefficients of a power series like \eqref{eqn:ps} with null radius of convergence. Define, for any $\gamma$ in $(0, 1]$,
\[
A_n(\gamma):=\max_{(1-\gamma)n\leq k\leq n}\alpha_k,
\]
and, for the sake of brevity,
\[ 
A_n:=\max_{1\leq k\leq n}\alpha_k.
\]
 Hence define 
\[
L_n(\gamma):=\frac{A_n(\gamma)}{A_n},
\]
that is less then or equal to $1$, for any natural $n$ and for every $\gamma$. Looking for a way to measure the concentration of big terms of $(a_n)_{n\in\N}$ for indices lower but not too far from $n$, we define
\[
L(\gamma):=\liminf_{n\to\infty}L_n(\gamma),
\]
that is, in particular, an increasing function of $\gamma$. That makes the limit
\[
G:=\lim_{\gamma\to0^+}L(\gamma)
\]
exist: we will denote it as the {\bf gauge} of the sequences $(a_n)_{n\in\N}$, or, depending on the context, of the power series $f(z)=\sum_{n=0}^\infty a_nz^n$. 
\begin{ex}
The non convergent power series
\begin{equation*}
\label{eqn:ex1}
f_\rho(z):=\sum_{n=0}^\infty \rho^n!z^{\rho^n}
\end{equation*}
 has, for any $\rho=2, 3, \dots$, null gauge. That is due to the fact that the \emph{gaps} in the coefficients of $f$ are big enough to make, for any positive $\gamma$, $A_n(\gamma)$ being zero for infinitely many $n$. As one can notice, that is the smallest magnitude of gaps we have to require in order to get $G=0$. For instance, a power series like 
 \begin{equation*}
 \label{eqn:ex2}
 f_k(z):=\sum_{n=0}^\infty n^k!z^{n^k}
 \end{equation*}
 has null radius of convergence and  gauge equal to $1$. Observe also that, fixed $t$ in $[0, 1]$, it is possible to find a power series with $R=0$ and gauge $G=t$. Just define the $n$-th coefficient of $f$ as
 \[
 a_n:=
 \begin{cases}
 \begin{split}
 n^n\qquad& \text{if}\quad n\in\set{2^k | k\in\N}\\
 1+(tn+\sqrt{n})^n\qquad&\text{otherwise}
 \end{split}
 \end{cases}.
 \]
 This also shows that a non convergent power series can have null gauge even with no zero coefficients.
 \end{ex}
A non convergent power series with gauge $G<1$ is the analogous of a power series with positive and finite radius of convergence that has Ostrowski gaps in its coefficients:
 \begin{defi}
\label{def:ostr}
Let $f=\sum_{n=0}^\infty a_nz^n$ be a power series with radius of convergence $R$. Then $f$ has \emph{Ostrowski gaps} if there exists a positive $\gamma$ such that
\begin{equation}
\label{ostr}
\liminf_{n\to\infty}A_n(\gamma)<\frac{1}{R}.
\end{equation}
\end{defi}
 As we will see in Section \ref{sec:ps}, the characterization of the Szegö class for non convergent power series will be given in terms of their gauge, while the well known result for $0<R<\infty$ (Theorem \ref{teo:cb}) is given in terms of the presence of Ostrowski gaps.

\subsection{Zeros of polynomials}
Given the nature of the problem we are studying, it is natural to look for relations between zeros of polynomials and their coefficients. In particular, we are looking for bounds for the absolute value of those zeros, since we are going to see if either they accumulate to some circle, or they don't. Here we report some classical results. It is worth to observe from now that those are bounds that apply to any polynomials and not only to sections of a prefixed power series. The result that it is going to be useful for the proof of one implication of Theorem \ref{theo:main}, Lemma \ref{lemma:vv}, can be found in \cite{josechu}.\\

Let $P(z)=b_0+\dots+b_nz^n$ be a polynomial of degree $n$. A classical result due to Cauchy asserts that all the zeros of $P$ lie in $\overline{\D}_{C_P}$, where $C_{P}$ is the unique positive root of 
\begin{equation}
\label{eqn:cauchy}
\abs{b_n}x^n=\sum_{k=0}^{n-1}\abs{b_k}x^k.
\end{equation}
A generalization of this result is given by the so called {\bf Van Vleck's bounds}: $P$ has actually at least $m$ zeros in $\overline{\D}_{V_{P, m}}$, where $V_{P, m}$ is the unique positive root of
\begin{equation}
\label{eqn:vv}
\abs{b_n}x^n=\sum_{j=0}^{m-1}\binom{n-j-1}{m-j-1} \abs{b_j}x^j.
\end{equation}
Observe that for $n=m$ equations \eqref{eqn:cauchy} and \eqref{eqn:vv} coincide. For a more detailed description of Cauchy's and Van Vleck's bounds, see \cite{zpoly}, chapters VII and VIII.\\
To get a lower bound for the position of some roots of $P$, recall that the \emph{reverse companion polynomial} of $P$ is defined by relocating backwards the coefficients of $P$, that is
\[
Q(z):=z^nP\bigg(\frac{1}{z}\bigg)=b_n+b_{n-1}z+\dots+b_0z^n.
\]
The zeros of $Q$ are the reciprocal of the zeros of $P$.\\
 It is also going to be useful the following well known upper bound for binomial coefficients: for any $n\geq 1$ and $0\leq k\leq n$
\begin{equation}
\label{eqn:binom}
\binom{n}{k}\leq e^{nH(k/n)},
\end{equation}
where $H(x)$ is the so called \emph{entropy function}
\[
H(x):=x\log{\frac{1}{x}}+(1-x)\log{\frac{1}{1-x}}, \qquad  0\leq x\leq1.
\]
Inequality in \eqref{eqn:binom} can be easily derived by a strong version of Stirling's formula proven by Robbins in \cite{rob}.
Observe that $H$ attain its maximum at $x=1/2$ and $H(0)=H(1)=0$.\\
 The idea that is going to be useful for our purpose is that given a polynomial $P$ and $m\leq \deg(P)$, we can find a radius $v$ so that at least $m$ zeros of $P$ have modulus greater than $v$, and the relation between $m$, $v$ and the coefficients of $P$ is the following:
\begin{lemma}
\label{lemma:vv}
For any $m\leq n$ and polynomial $P(z)=\sum_{k=0}^nb_nz^n$ such that $b_0\ne0$, there exists a positive radius $v_{P, m}$ so that
\begin{equation}
\label{eqn:pest}
\abs{b_0}\leq e^{nH((m-1)/n)}\max_{n-m+1\leq k\leq n}\abs{b_k}\cdot\max\set{1, v_{P, m}}^n
\end{equation}
and $P$ has at least $m$ zeros whose absolute values exceed $v_{P, m}$.
\begin{proof}
Let $Q$ be the reverse companion polynomial of $P$. Then Van Vleck's bounds say that $Q$ has at least $m$ zeros in $\overline{\D}_\alpha$, where $\alpha$ is the unique positive root of 
\begin{equation}
\label{eqn:vvQ}
\abs{b_0}x^n=\sum_{j=0}^{m-1}\binom{n-j-1}{m-j-1}\abs{b_{n-j}}x^j.
\end{equation}
Since the zeros of $P$ are the reciprocal of the ones of $Q$, $P$ has at least $m$ zeros outside the open disc of radius $v_{P, m}:=\alpha^{-1}$, and by uniqueness of $\alpha$ as a positive root of \eqref{eqn:vvQ}, $v_{P, m}$ is the unique positive root of the reversed companion polynomial of \eqref{eqn:vvQ}, that is
\[
\abs{b_0}=\sum_{k+n-m+1}^n\binom{k-1}{k-(n-m)-1}\abs{b_k}x^k.
\]
In particular then
\[
\begin{split}
\abs{b_0}&\leq\big(\max_{n-m+1\leq k\leq n}\abs{b_k}\big)\max\set{1, v_{P, m}}^n\sum_{k=n-m+1}^n\binom{k-1}{k-(n-m)-1}\\
&=\binom{n}{m-1}\big(\max_{n-m+1\leq k\leq n}\abs{b_k}\big)\max\set{1, v_{P, m}}^n\\
&\leq e^{nH((m-1)/n)}\big(\max_{n-m+1\leq k\leq n}\abs{b_k}\big)\max\set{1, v_{P, m}}^n,
\end{split}
\]
as we claimed.
\end{proof}
\end{lemma}
One can observe that the maximum involving the coefficients in \eqref{eqn:pest} considers only some coefficients before the $n$th, and that we can choose the number $m$ of this coefficients to be as small as we want. It is, therefore, a very suitable estimation for the notion of gauge we saw in section \ref{sec:gauge}. 
\subsection{A Jentzsch-Szegö Type Theorem}
Let $E$ be a compact subset of $\C$ such that $\Omega:=\overline{\C}\setminus E$ is connected. Then it is well known that that there exists a Green function of $\Omega$ with pole at $\infty_{\C}$, i.e., a non negative harmonic  function $g$ such that
\[
\lim_{z\to\de E}g(z)=0
\]
and
\[
\lim_{z\to\infty}(g(z)-\log\abs{z})=-\log\cp(E).
\]
Here, by $\cp(E)$ we mean the logarithmic capacity of $E$, defined by
\[
\cp(E):=e^{-V(E)},
\]
where $V(E)$ is the infimum of the energy quantity
\[
\int_\C\int_\C\log\frac{1}{\abs{w-z}}\, d\sigma(z)\,d\sigma(w)
\]
along all the positive Borel measures $\sigma$ supported on $E$. We will assume that $\cp(E)$ is strictly positive, since in this case there exists a unique measure $\mu_E$ supported on $E$ such that
\[
V(E)=\int_\C\int_\C\log\frac{1}{\abs{w-z}}\, d\mu_E(z)\,d\mu_E(w).
\]
The measure $\mu_E$ is called \emph{equilibrium measure} of $E$. For instance, see \cite[p. 10, 11]{andri}
\begin{ex}
\label{ex:disc}
The logarithmic capacity $\cp(\overline{\D})$ of the closure of the unit disc is $1$, since the Green function in this case is
\[
g(z)=\log\abs{z},
\]
 and its equilibrium measure is $\Lambda$. For a proof, see \cite[p. 12]{andri}.
\end{ex}
The following result, whose proof can be found in \cite[p. 50-55]{andri}, looks at the asymptotic distribution of the zeros of an arbitrary sequence of polynomials at the boundary of a compact and simply connected subset of the complex plane:
\begin{theo}
\label{teo:sj}
Let $E$, $\Omega$ and $g$ be as above, and let $(p_n)_{n\in\N}$ be a sequence of complex valued polynomials, with associated sequence of zeros counting measures $(\nu_n)_{n\in\N}$ . Suppose that there exists a sequence of integers $(k_n)_{n\in\N}$ such that $k_n\geq\deg(p_n)$ and 
\begin{align*}
&\limsup_{n\to\infty}\frac{1}{k_n}\max_{z\in E}\log\abs{p_n(z)}\leq0\\
&\liminf_{n\to\infty}\bigg[\max_{z\in S}\bigg(\frac{1}{k_n}\log\abs{p_n(z)}-g(z) \bigg) \bigg]\geq0,
\end{align*}
for some compact subset $S$ of $\Omega$. Then:
\begin{description}
\item[(i)] For any compact subset $K$ of $\Omega$
\[
\lim_{n\to\infty}\nu_n(K)=0;
\]
\item[(ii)] If for any compact subset $M$ of the interior of $E$ 
\[
\lim_{n\to\infty}\nu_n(M)=0,
\]
then $\nu_n$ converges weakly to $\mu_E$.
\end{description}
\end{theo}
\begin{obs}
This implies, for example, that the equilibrium measure of $E$ is supported on its boundary.
\end{obs}
For the purposes of this work, $E$ is going to be the closure of the unit disk. In this case, thanks to Example \ref{ex:disc}, Theorem \ref{teo:sj} can be restated as follows:
\begin{coro}
\label{coro:sj}
 Let $(p_n)_{n\in\N}$ be a sequence of complex valued polynomials, with associated sequence of zeros counting measures $(\nu_n)_{n\in\N}$ . Suppose that there exists a sequence of integers $(k_n)_{n\in\N}$ such that $k_n\geq\deg(p_n)$ and 
\begin{align}
\label{eqn:teo:a}
&\limsup_{n\to\infty}\frac{1}{k_n}\max_{z\in \overline{\D}}\log\abs{p_n(z)}\leq0\\
\label{eqn:teo:b}
&\liminf_{n\to\infty}\max_{z\in S}\frac{\abs{p_n(z)}^\frac{1}{k_n}}{\abs{z}}  \geq1,
\end{align}
for some compact subset $S$ of $\set{\abs{z}>1}$. Then:
\begin{description}
\item[(i)] For any compact subset $K$ of $\overline{\C}\setminus\overline{\D}$
\[
\lim_{n\to\infty}\nu_n(K)=0;
\]
\item[(ii)] If for any compact subset $M$ of $\D$ 
\[
\lim_{n\to\infty}\nu_n(M)=0,
\]
then $\nu_n$ converges weakly to $\Lambda$.
\end{description}
\end{coro}
\section{Zeros of sections of power series}
\label{sec:ps}
Suppose that the radius of convergence $R$ of 
\[
f(z)=\sum_{n=0}^\infty a_nz^n 
\]
is positive and finite. By eventually replacing $f$ with $z\mapsto f(z/R)$, we can assume that $R=1$. Let $(\nu_n)_{n\in\N}$ be the zero counting measure of $f$, and
\[
\abs{\nu_n}:=\frac{1}{n}\sum_{w\in Z_n}\delta_{\abs{w}}.
\]
We say that $(\nu_n)_{n\in\N}$ \emph{converges absolutely} to $\de\D$ if $(\abs{\nu_n})_{n\in\N}$ converges weakly to $\delta_1$.
Thanks to Hurwitz's Theorem, for any $r<1$
\begin{equation}
\label{hur}
\lim_{n\to\infty}\abs{\nu_n}([0, r])=0.
\end{equation}
Without assuming any regularity for the coefficients, the weak convergence of $(\nu_n)_{n\in\N}$ to $\Lambda$ holds only for a subsequence:
\begin{theo}
\label{theo:js}
Let $f$ be a power series with radius of convergence $R=1$, and let $(\nu_n)_{n\in\N}$ be its sequence of zero counting measures. Then:
\begin{description}
\item[(i)] there exists a subsequence $(\nu_{n_k})_{k\in\N}$ that converges absolutely to $\de\D$;
\item[(ii)] if a subsequence  $(\nu_{n_k})_{k\in\N}$ converges absolutely to $\de\D$, then it converges weakly to $\Lambda$.
\end{description}
\end{theo}
Part (i) is due to Jentszch \cite{jen}, while the outstanding improvement in (ii) was proven by Szegö in \cite{szego}: the fact that the zeros of the partial sums of $f$ approach asymptotically the unit circle, implies that they do that in an equidistributed way. \\
For the result to be true for the whole sequence of zero counting measure, it turns out that a necessary and sufficient condition is that the power series has not big gaps in its coefficients. The following characterization of the Szegö class was first announced by Carlson in \cite{carlo}, while a complete proof appeared the first time in the monograph \cite{bou} of Bourion. An alternative proof is given in \cite{erdos} by Erdös and Fried, who also credited the following Theorem to Bourion:
\begin{theo}[Carlson-Bourion]
\label{teo:cb}
If $R=1$, then $f$ belongs to $\s$ if and only if it has not Ostrowski gaps.
\end{theo}
In particular, Erdös and Fried proved that $f$ has Ostrowsi gaps if and only if there exists $\rho>1$ so that 
\[
\liminf_{n\to\infty}\,\, \abs{\nu_n}([0, \rho]])<1
\]
which, together with \eqref{hur}, says that the zero counting measure of $f$ converges absolutely to $\de\D$ if and only if $f$ has no Ostrowsi gaps. Thanks to Theorem \ref{theo:js}, part (ii), this is equivalent to say that $f$ belongs to $\s$.

\subsection{The case $R=0$}
Suppose now that $R=0$, i.e., that 
\begin{equation}
\label{eqn:r0}
\limsup_{k\to\infty}\alpha_k=+\infty.
\end{equation}
The idea  is to normalize the sections $s_n$, and to look for an equidistribution of their zeros along $\de\D$. More precisely, in order to give an analogous of  Theorem \ref{teo:cb} for the case $R=0$, we need to find a sequence $(r_n)_{n\in\N}$ of positive numbers such that the sequence of polynomials defined by
\[
P_n(z)=s_n\bigg(\frac{z}{r_n}\bigg)
\]
belongs to $\s$. The reason why we need to normalize, at each step $n$, the $n$th partial sum of $f$ can be found  in \cite[Theorem 1]{rubel}, where Dilcher and Rubel showed that, for any positive $\varepsilon$ and for $n$ large enough, the zeros of the sequence of polynomials
\begin{equation*}
s_{k_n}(f)\bigg(\frac{z}{\alpha_{k_n}}\bigg)
\end{equation*}
are all contained in the annulus
\[
\set{\abs{a_0}/(1+\abs{a_0})\leq\abs{z}\leq2+\varepsilon},
\]
where $(k_n)_{n\in\N}$ is the infinite sequence such that
\[
\alpha_{k_n}=A_{k_n}.
\]
This, in particular, implies that $(\nu_{k_n})_{n\in\N}$ converges to the Dircac delta at the origin. They also showed that zeros of the normalized sequence of the sections
\begin{equation}
\label{eqn:norm1}
s_n\bigg(\frac{z}{\alpha_n}\bigg)
\end{equation}
converges absolutely to $\de\D$ under some conditions on the coefficients:
\begin{theo}[Dilcher, Rubel]
\label{teo:dr}
Let $f$ be a power series like \eqref{eqn:ps} with null radius of convergence and $a_0\ne0$. Suppose that 
\begin{enumerate}
\item $a_n$ is real and non negative for any integer $n$,
\item $a_{n-1}a_{n+1}>a_n^2$ for all $n\geq1$.
\end{enumerate}
Then the zeros counting measures of the normalized sections in \eqref{eqn:norm1} converge absolutely to $\de\D$.
\end{theo}
One can observe that the result does not necessary imply that the roots of the normalized sections are asymptotically equidistributed on $\de\D$, since Theorem \ref{theo:js}, part (ii) holds only for sections of a predetermined power series with positive and finite radius of convergence.\\

Here we generalize Theorem \ref{teo:cb} at the case $R=0$, by considering a normalization factor slightly different from \eqref{eqn:norm1}. Indeed, if one considers the normalized sequence
\begin{equation}
\label{eqn:norm2}
p_n(z):=s_n\bigg(\frac{z}{A_n} \bigg)=\sum_{k=0}^n a_kA_n^{-k}z^k
\end{equation}
then the following restatement of Theorem \ref{theo:main} in terms of the gauge of $f$ holds:
\begin{theo}
\label{teo}
Let
\[
f(z)=\sum_{n=0}^\infty a_nz^n
\]
be a power series with null radius of convergence. Then the sequence of the normalized section of $f$ in \eqref{eqn:norm2} belongs to the Szegö class if and only if the gauge of $f$ is $1$.
\end{theo}
Observe that
\begin{itemize}
\item even if the normalization factor is different from the one in Theorem \ref{teo:dr}, we are supposing only that the gauge of $f$ is $1$, that is a more general condition than properties 1 and 2 of  Theorem \ref{teo:dr}. Indeed, condition 2 implies that the coefficients are not null and
\[
a_{n+1}>\frac{a_n}{a_{n-1}}a_n>\frac{a_{n-1}}{a_{n-2}}a_{n}>\dots>\frac{a_1}{a_0}a_n,
\]
that implies that the gauge of $f$ is $1$.
\item the normalization factor is not so much different for the one in Theorem \ref{teo:dr}: as we have already observed, there exists an infinite sequence $(k_n)_{n\in\N}$ such that
\[
\alpha_{k_n}=A_{k_n},
\]
and this thank to \eqref{eqn:r0}. Thus, looking at the normalized section in \eqref{eqn:norm1}, Theorem \ref{teo} implies that there exists a subsequence of the normalized zero counting measures that converges at $\Lambda_1$, while, as we have already pointed out, Theorem \ref{teo:dr} does not say nothing about the equidistribution of the zeros of the sequence in \eqref{eqn:norm1}.
\end{itemize}
\begin{proof}
Suppose that $G=1$. We want to apply Corollary \ref{coro:sj} at the sequence $(p_n)_{n\in\N}$, with $k_n=n$. Observe that, for any natural $1\leq k\leq n$ and positive $\delta$
\[
\bigg(\frac{\alpha_k}{A_n} \bigg)^k=\frac{\abs{a_k}}{A_n^k}=\babs{\frac{1}{2\pi i}\int_{\abs{z}=1+\delta}\frac{p_n(z)}{z^{k+1}}\, dz}\leq\max_{\abs{z}=1+\delta}\frac{\abs{p_n(z)}}{\abs{z}^k}, 
\]
so that, for all $1\leq k\leq n$,
\[
\bigg(\frac{\alpha_k}{A_n}\bigg)^\frac{k}{n}\leq\max_{\abs{z}=1+\delta}\frac{\abs{p_n(z)}^\frac{1}{n}}{\abs{z}^\frac{k}{n}}.
\]
Then, for any positive $\gamma$, 
\[
\begin{split}
&\max_{\abs{z}=1+\delta}\frac{\abs{p_n(z)}^\frac{1}{n}}{\abs{z}}\geq\max_{k=1,\dots,n}(1+\delta)^{\frac{k}{n}-1}\bigg( \frac{\alpha_k}{A_n}\bigg)^\frac{k}{n}\\
\geq&\max_{(1-\gamma)n\leq k\leq n}(1+\delta)^{\frac{k}{n}-1}\bigg( \frac{\alpha_k}{A_n}\bigg)^\frac{k}{n}\geq(1+\delta)^{-\gamma}L_n(\gamma).
\end{split}
\]
Thus
\[
\liminf_{n\to\infty}\max_{\abs{z}=1+\delta}\frac{\abs{p_n(z)}^\frac{1}{n}}{\abs{z}}\geq(1+\delta)^{-\gamma}L(\gamma),
\]
and letting $\gamma$ going to $0$ we get \eqref{eqn:teo:b}, since $G=1$.\\
To see that \eqref{eqn:teo:a} holds, it suffices to observe that 
\begin{equation}
\label{eqn:max}
\max_{z\in\overline{\D}}\abs{p_n(z)}=\max_{z\in\de\D}\abs{p_n(z)}\leq\sum_{k=0}^n\abs{a_k}A_n^{-k}\leq\abs{ a_0}+n\cdot\max_{k=1, \dots, n}\bigg(\frac{\alpha_k}{A_n} \bigg)^k=\abs{a_0}+n.
\end{equation}
To conclude, we have to show that that, for any compact subset $M$ of $\D$,
\begin{equation}
\label{eqn:0ins}
\lim_{n\to\infty}\nu_n(M)=0.
\end{equation}
Unfortunately, we can not use Hurwitz's Theorem, since $(p_n)_{n\in\N}$ is not supposed to converge uniformly on compact subsets of $\D$. We may notice that it suffices to show that \eqref{eqn:0ins} holds for any disc $\overline{\D}_r$ centered at zero and of radius $r$. Let then $r<1$, and observe that
\[
\frac{1}{n}\sum_{w\in Z_{p_n}\cap\overline{\D}r}\log\frac{1}{\abs{w}}\geq\bigg(\log\frac{1}{r}\bigg)\nu_n(\overline{\D_r}).
\]
Also, without loss of generality, we can assume that $a_0\ne0$. Indeed, if $k$ is the multiplicity of $z=0$ as a zero of $f$, one can replace $f(z)$ with $f(z)/z^k$, and an eventual weak limit $\nu$ of the zero counting measures wouldn't be affected, since we are removing only finitely many zeros. Thus Jensen's formula and the inequality above yield
\[
\begin{split}
\nu_n(\overline{\D}_r)\log{\frac{1}{r}}&\leq\frac{1}{n}\sum_{w\in Z_{p_n}\cap\D}\log{\frac{1}{\abs{w}}}=\frac{1}{2\pi}\int_0^{2\pi}\frac{1}{n}\log{\frac{\abs{p_n(e^{i\theta})}}{\abs{a_0}}}\, d\theta\\
&\leq\frac{1}{n}\max_{z\in\de\D}\log{\frac{\abs{p_n(z)}}{\abs{a_0}}}\leq \frac{1}{n}\log{\frac{\abs{a_0}+n}{\abs{a_0}}},
\end{split}
\]
where the last inequality holds thanks to \eqref{eqn:max}. Therefore,
\[
\lim_{n\to\infty}\nu_n(\overline{\D}_r)=0,
\]
for any $r<1$, as claimed.\\

Conversely, suppose that $G<1$, fix \footnote{we may choose any $1<r<\infty$ if $G$ is null}$1<r<1/G$, and a positive $\varepsilon$ so that $r(G+\varepsilon)<1$. Since the entropy function $H$ is continuous and null at $x=0$, we can choose a $\gamma$ in $(0, 1/2)$  so that
\[
L(\gamma)^{1-\gamma}e^{H(\gamma)}\leq G+\varepsilon.
\]
Thanks to Lemma \ref{lemma:vv}, each $p_n$ has at least $m_n:=\floor{\gamma n}+1$ zero outside the disc with center $z=0$ and radius $v_n$, where 
\[
\begin{split}
\abs{a_0}&\leq e^{nH(\floor{\gamma n}/n)}\max_{(1-\gamma)n\leq k\leq n}\bigg(\frac{\alpha_k}{A_n} \bigg)^k\max\set{1, v_n}^n\\
&\leq  e^{nH(\floor{\gamma n}/n)}L_n(\gamma)^{n(1-\gamma)}\max\set{1, v_n}^n,
\end{split}
\]
that yields
\[
\limsup_{n\to\infty}\big(\max\set{1, v_n}\big)\geq\frac{1}{e^{H(\gamma)}L(\gamma)^{1-\gamma}}>\frac{1}{G+\varepsilon}>r,
\]
thanks to how we have chosen $\varepsilon$ and $r$. Since $r>1$, this implies that, for infinitely many $n$, $v_n$ is greater than $r$ and $p_n$ has at most $n-\floor{\gamma n}-1$ zeros in $\D_{r}$. Therefore
\[
\liminf_{n\to\infty}\nu_n(\de\D)\leq\liminf_{n\to\infty}\nu_n(\D_r)\leq1-\gamma<1,
\]
and so $(p_n)_{n\in\N}$ does not belong to the Szegö class.
\end{proof}
The argument for the second half of the proof can be also found in \cite{josechu}.\\

One can observe that we never used, in the proof of Theorem \ref{teo}, the fact that $R=0$. Therefore, the same argument holds for a power series with positive and finite radius of convergence. The fact is that, in this case, the gauge of $f$ being $1$ doesn't translate to $f$ not having  Ostrowski gaps, and therefore Theorem \ref{teo} is not equivalent to Theorem \ref{teo:cb}. Of course, both those conditions are measuring the density of the big terms in the sequence of the coefficients of $f$. Nevertheless, the way in which one wants to look at such a condition depends on the finiteness of 
\[
\limsup_{n\to\infty}\alpha_n,
\]
and so on the fact that the radius of convergence of $f$ is null or positive. Theorem \ref{teo} shows that, once we choose the normalization factor in \eqref{eqn:norm2}, the right way to measure Ostrowski-type gaps for a power series with null radius of convergence is the notion of gauge defined in Section \ref{sec:gauge}, if we want to give an analogous of Theorem \ref{teo:cb}. A notion of gauge that translates to the presence of Ostrowski gaps for the case $R=1$  can be found in \cite{josechu}, together with a characterization of the Szegö class in terms of such a gauge.


\begin{thebibliography}{99}
\bibitem{andri} \textsc{Andrievskii, V. V. and Blatt, H. P.:} \emph{Discrepancy of Signed Measures and Polynomial Approximation}, Springer-Verlag, New York, 2002.
\bibitem{bou} \textsc{Bourion, G.:} L'ultraconvergence dans les séries de Taylor, \emph{Actualités scientifiques et industrielles}, no. 472, Paris, 1937.
\bibitem{carlo} \textsc{Carlson, F.:} Sur quelques suites de polynomes, \emph{C. R. Acad. Sci. Paris} 178 (1924), 1677-1680
\bibitem{rubel} \textsc{Dilcher, K. and Rubel, L. A.:} Zero Section of Divergent Power Series, \emph{J. Math. Anal. Appl.} 198 (1996), no.1, 98-110.
\bibitem{erdos} \textsc{Erdös, P. and Fried, H.:} On the connection between Gaps in Power series and the Roots of Their Partial Sums, \emph{Trans. of the American Mathematical Society}, vol. 62 no.1, 1947, 53-61.
\bibitem{josechu} \textsc{Fern\'andez, J. L.:} Zeros of Sections of Power Series: Deterministic and Random, \emph{Computational Methods and Function Theory}, September 2017, Volume 17, Issue 3, 463-486
\bibitem{jen} \textsc{Jentzsch, R.:} Untersuchungen zur Theorie der Folgen analytischer Funktionen, \emph{Acta Math.} 41 (1918), 219
\bibitem{zpoly} \textsc{Marden, M.:} \emph{Geometry of polynomials}, Mathematical Surveys and Monographs, no 3, American Mathematical Society, 1966.
\bibitem{rob} \textsc{Robbins, H.:} A remark on Stirling's Formula, \emph{The American Mathematical Monthly}, Vol 62, no 1, 1955, 26-29.
\bibitem{szego} \textsc{Szegö, G.:} Über die Nullstellen von Polynomen, die in einem Kreis gleichmässig konvergieren, \emph{Sitzungsber. Berliner Math. Ges.}, 21 (1922), 59-64
\end{thebibliography}
\end{document}